\documentclass[12pt]{iopart}
\usepackage{amssymb}
\usepackage{amsthm}
\usepackage{graphicx}
\usepackage{graphics,epsfig,color}

\newcommand{\bfx}{{\bf x}}

\newcommand{\N}{{\mathbf N}}

\newcommand{\R}{{\mathbf R}}

\newcommand{\Z}{{\mathbf Z}}
\newcommand{\C}{{\mathbf C}}

\newcommand{\bfxz}{{\mathbf{x}}_0}
\newcommand{\sign}{{\rm sign}}
\newcommand{\ce}{{\rm ce}}
\newcommand{\se}{{\rm se}}

\newcommand{\me}{{\rm me}}

\newcommand{\Ie}{{\rm Ie}}
\newcommand{\Io}{{\rm Io}}
\newcommand{\Ke}{{\rm Ke}}
\newcommand{\Ko}{{\rm Ko}}

\newtheorem{lemma}{Lemma}[section]
\newtheorem{thm}[lemma]{Theorem}

\begin{document}
\bibliographystyle{plain}

\title[Eigenfunction expansions in parabolic and elliptic cylinder 
coordinates]{Eigenfunction expansions for a fundamental solution of Laplace's equation
on $\R^3$ in parabolic and elliptic cylinder coordinates}

\author{H S Cohl${}^{1}$ and H Volkmer${}^2$}

\address{$^1$Information Technology Laboratory, National Institute of Standards and Technology, Gaithersburg, MD, USA}
\address{$^2$Department of Mathematical Sciences, University of Wisconsin--Milwaukee, 
P.~O. Box 413, Milwaukee, WI 53201, U.S.A.}
\ead{hcohl@nist.gov}
\begin{abstract}
A fundamental solution of Laplace's equation in three dimensions is expanded
in harmonic functions that are separated in parabolic or elliptic cylinder coordinates.
There are two expansions in each case which reduce to expansions of the Bessel 
functions $J_0(kr)$ or $K_0(kr)$, $r^2=(x-x_0)^2+(y-y_0)^2$, in parabolic and elliptic 
cylinder harmonics. Advantage is taken of the fact that $K_0(kr)$ is a fundamental solution
and $J_0(kr)$ is the Riemann function of partial differential equations on the Euclidean plane.
\end{abstract}

\pacs{02.30.Em, 02.30.Gp, 02.30.Hq, 02.30.Jr, 02.30.Mv, 02.40.Dr}
\ams{35A08, 35J05, 42C15, 33E10, 33C15}
\maketitle

\section{Introduction}
A fundamental solution of Laplace's equation
\begin{equation}
\frac{\partial^2 U}{\partial x^2}+\frac{\partial^2 U}{\partial y^2}+\frac{\partial^2 U}{\partial z^2}=0
\label{1:Laplace}
\end{equation}
is given by (apart from a multiplicative factor of $4\pi$)
\begin{equation}
U(\bfx,\bfxz)=\frac{1}{\|\bfx-\bfxz\|}, \quad \textup{where
$\bfx=(x,y,z)\ne \bfxz=(x_0,y_0,z_0)$} , \label{1:fund}
\end{equation}
and $\|\bfx-\bfx_0\|$ denotes the Euclidean distance between $\bfx$ and $\bfx_0$.
In many applications it is required to expand a fundamental solution in the form of a series 
or an integral, in terms of
solutions of (\ref{1:Laplace}) that are separated in suitable curvilinear coordinates.
Examples of such applications include electrostatics, magnetostatics, quantum direct and exchange Coulomb 
interactions, Newtonian gravity, potential flow, and steady state heat transfer.
Morse \& Feshbach (1953) \cite{MorseFesh}
(see also Hobson (1955) \cite{Hob}; MacRobert (1947) \cite{MacRobert47}; Heine (1881) \cite{Heine})
provide a list of such expansions for various coordinate systems but the formulas
for several coordinate systems are missing. It is the goal of this paper to provide these expansions for
parabolic and elliptic cylinder coordinates.
Although these expansions are partially known from
Buchholz (1953) \cite{Buchholz}
and Lebedev (1972) \cite{Lebedev} for parabolic cylinder coordinates,
and from Meixner \& Sch{\"a}fke (1954) \cite{MeixnerSchafke54} for
elliptic cylinder coordinates,
we found it desirable to investigate these expansions in
a systematic fashion and provide direct proofs for them based
on eigenfunction expansions.

There will be two expansions for both of these
coordinate systems. The first expansion for a cylindrical coordinate
system on $\R^3$ starts from
the known formula in terms of the
integral of Lipschitz (see Watson (1944) 
\cite[section 13.2]{Watson}; Cohl {\it et al.} (2000) \cite[(8)]{CTRS})
\begin{equation}
\frac{1}{\|\mathbf{x}-\mathbf{x_0}\|}= \int_0^\infty
J_0(k \sqrt{(x-x_0)^2+(y-y_0)^2})e^{-k|z-z_0|}\,dk.
\label{1:firstexpansion}
\end{equation}
The Bessel function $J_\nu$ can be defined by
(see for instance (10.2.2) in Olver {\it et al.} (2010) \cite{NIST})
\begin{equation}
J_\nu(z):=\biggl(\frac{z}{2}\biggr)^\nu\sum_{n=0}^\infty \frac{(-z^2/4)^n}{n!\Gamma(\nu+n+1)},
\label{1:J}
\end{equation}
Note that for $W_k:\R^2\times\R^2\to\R$ 
\begin{equation}
W_k(x,y,x_0,y_0):=J_0(kr),
\label{1:h1}
\end{equation}
where $r^2:=(x-x_0)^2+(y-y_0)^2,$
solves the partial differential equation
\begin{equation}
\frac{\partial^2 U}{\partial x^2}+\frac{\partial^2 U}
{\partial y^2} +k^2 U =0.
\label{1:pde1}
\end{equation}
Therefore, in a cylindrical coordinate system on $\R^3$ involving
the Cartesian coordinate $z$, the first expansion
(\ref{1:firstexpansion}) reduces to
expanding $W_k$ from (\ref{1:h1}) in terms of solutions of (\ref{1:pde1}) that are separated in a curvilinear coordinate system on the plane.

The second expansion for a cylindrical coordinate system on $\R^3$ is based on
the known formula given in terms of the Lipschitz-Hankel integral
(see Watson (1944) \cite[section 13.21]{Watson}; Cohl {\it et al.} (2000) \cite[(9)]{CTRS})
\begin{equation}
\frac{1}{\|\mathbf{x}-\mathbf{x_0}\|}= \frac{2}{\pi} \int_0^\infty
K_0(k \sqrt{(x-x_0)^2+(y-y_0)^2})\cos k(z-z_0)\,dk ,
\label{1:secondexpansion}
\end{equation}
where $K_\nu:(0,\infty)\to\R$ (cf.~(10.32.9) in Olver {\it et al.} (2010)
\cite{NIST}), the modified Bessel function of the second kind
(Macdonald's function), of order $\nu\in\R$, is defined by
\begin{equation}\label{1:K}
K_\nu(z):=\int_0^\infty e^{-z\cosh t}\cosh(\nu t)\,dt.
\end{equation}
Now $V_k:\R^2\times\R^2\setminus\{(\bfx,\bfx):\bfx\in\R^2\}\to(0,\infty),$
defined by
\begin{equation}
V_k(x,y,x_0,y_0):=K_0(kr)
\label{1:h2}
\end{equation}
solves the partial differential equation
\begin{equation}
\frac{\partial^2 U}{\partial x^2}
+\frac{\partial^2 U}{\partial y^2} -k^2 U =0 .
\label{1:pde2}
\end{equation}
In a cylindrical coordinate system on $\R^3$ involving the Cartesian
coordinate $z$, the second expansion
(\ref{1:secondexpansion}) reduces to
expanding $V_k$ in terms of solutions of (\ref{1:pde2}) that are separated in curvilinear coordinates on the plane.

The paper is organized as follows. In section \ref{2} we derive
the desired expansion of $K_0(kr)$ in parabolic cylinder
coordinates. This expansion is given in terms of series over
Hermite functions. In section \ref{3} we obtain an integral
representation for $J_0(kr)$ in terms of separated solutions of
(\ref{1:pde2}) in terms of (modified) parabolic cylinder
functions. We show how these results are based on a general
expansion theorem in terms of solutions of the
differential equation
\begin{equation}\label{1:ode}
-u''-\frac14 \xi^2 u=\lambda u .
\end{equation}
In sections \ref{Jzeroellipticcylinder}
and \ref{ExpansionKzeroellipticcylinder}
we derive the fundamental
solution expansions in elliptic cylinder coordinates for
$J_0(kr)$ and $K_0(kr)$, respectively.

\medskip

Throughout this paper we rely on the following definitions.
The set of natural numbers is given by $\N:=\{1,2,3,\ldots\}$, the set
$\N_0:=\{0,1,2,\ldots\}=\N\cup\{0\}$, and the set
$\Z:=\{0,\pm 1,\pm 2,\ldots\}.$  The set $\R$ represents the real numbers
and the set $\C$ represents the complex numbers.

\section{Expansion of $K_0(kr)$ for parabolic cylinder coordinates}
\label{2}

Parabolic coordinates on the plane $(\xi,\eta)$ (see for
instance Chapter 10 in Lebedev (1972) \cite{Lebedev})
are connected to Cartesian coordinates $(x,y)$ by
\begin{equation}
x=\frac12(\xi^2-\eta^2),\quad y=\xi\eta,
\label{2:pc}
\end{equation}
where $\xi\in\R$ and $\eta\in[0,\infty)$.
To simplify notation we will first set $k=1$ in (\ref{1:h2}).
Then $V_1$ satisfies
\begin{equation}
\frac{\partial^2 U}{\partial x^2}+\frac{\partial^2 U}{\partial y^2}-U = 0
\quad \textup{if $(x,y)\ne (x_0,y_0)$.}
\label{2:pde1}
\end{equation}
Let $(\xi,\eta)$, $(\xi_0,\eta_0)$ be parabolic coordinates
on $\R^2$ for $(x,y)$ and $(x_0,y_0)$, respectively.
Then $V_1$ transforms to
\[
v(\xi,\eta,\xi_0,\eta_0):=K_0(r(\xi,\eta,\xi_0,\eta_0)),\quad (\xi,\eta)\ne \pm(\xi_0,\eta_0)
\]
and $r(\xi,\eta,\xi_0,\eta_0)$ is defined by
\begin{equation}\label{2:r}
r^2=\frac14\left[(\xi+\xi_0)^2+(\eta+\eta_0)^2\right]
\left[(\xi-\xi_0)^2+(\eta-\eta_0)^2\right],\quad r>0.
\end{equation}
Here and in the following we allow all $\xi,\eta,\xi_0,\eta_0\in\R$ such that $(\xi,\eta)\ne \pm(\xi_0,\eta_0)$.
From (\ref{2:pde1}) or by direct computation we obtain that $v$ solves the equation
\begin{equation}
\frac{\partial^2 u}{\partial \xi^2} +\frac{\partial^2 u}
{\partial \eta^2}- (\xi^2+\eta^2) u = 0.
\label{2:pde3}
\end{equation}
Separating variables $u(\xi,\eta)=u_1(\xi)u_2(\eta)$ in (\ref{2:pde3}),
we obtain the ordinary differential equations
\begin{eqnarray}
\label{2:ode1}
u_1''+(2n+1-\xi^2) u_1 & = & 0,\\
u_2''-(2n+1+\eta^2) u_2 & = & 0,
\label{2:ode2}
\end{eqnarray}
where we will use only $n\in\N_0$.
Equations (\ref{2:ode1}), (\ref{2:ode2}) have the general solutions
\begin{eqnarray*}
 u_1(\xi)& = & c_1 e^{-\xi^2/2} H_n(\xi) +c_2 e^{\xi^2/2}H_{-n-1}(i\xi),\\
 u_2(\eta)& =& c_3 e^{\eta^2/2} H_n(i\eta)+c_4e^{-\eta^2/2} H_{-n-1}(\eta) ,
\end{eqnarray*}
where $H_\nu:\C\to\C$ is the Hermite function which can be defined
in terms of Kummer's function of the first kind
$M$ as
(cf.~(10.2.8) in Lebedev (1972) \cite{Lebedev})
\[
H_\nu(z):=
\frac{2^\nu\sqrt{\pi}}{\Gamma\left(\frac{1-\nu}{2}\right)}
M\left(-\frac{\nu}{2},\frac12,z^2\right)
-
\frac{2^{\nu+1}\sqrt{\pi}}{\Gamma\left(-\frac{\nu}{2}\right)}
z M\left(\frac{1-\nu}{2},\frac32,z^2\right),
\]
and
\begin{equation}
M(a,b,z):=\sum_{n=0}^\infty\frac{(a)_nz^n}{(b)_nn!}
\label{Kummerfirst}
\end{equation}
(see for instance (13.2.2) in Olver {\it et al.} (2010) \cite{NIST}).
Note that the Kummer function of the first kind is entire
in $z$ and $a,$ and is a meromorphic function of $b.$
The Hermite function is an entire function of both $z$ and $\nu$.
If $\nu=n\in\N_0$ then $H_\nu(z)$ reduces to the Hermite polynomial
of degree $n$.

We will expand $v$ as a function of $\xi$ in an orthogonal series of
functions $e^{-\xi^2/2} H_n(\xi)$, $n\in\N_0$, so that the coefficients
\[
f_n(\eta,\xi_0,\eta_0):=\int_{-\infty}^\infty v(\xi,\eta,\xi_0,\eta_0) e^{-\xi^2/2} H_n(\xi)\,d\xi ,
\]
are to be evaluated.  We do this based on the observation that
the function $(\xi,\eta)\mapsto v(\xi,\eta,\xi_0,\eta_0)$ is a fundamental solution of
equation (\ref{2:pde3}). It has logarithmic singularities at the points $\pm(\xi_0,\eta_0)$.
Arguing as in Volkmer (1984) \cite[Theorem 1.11]{Volkmer84},
we obtain the following integral representation for solutions of (\ref{2:pde3}).

\begin{thm}
Let $u\in C^2(\R^2)$ be a solution of (\ref{2:pde3}). Let $(\xi_0,\eta_0)\in\R^2$,
and let $C$ be a a closed rectifiable curve on $\R^2$ which does not pass through
$\pm (\xi_0,\eta_0)$, and let $n^\pm$ be the winding number of $C$ with respect
to $\pm (\xi_0,\eta_0)$. Then we have
\[
\fl 2\pi \left[n^+u(\xi_0,\eta_0)+n^-u(-\xi_0,-\eta_0)\right]
=\int_C (u\partial_2 v-v\partial_2 u)\,d\xi+(v\partial_1 u-u\partial_1 v)\,d\eta,
\]
where $\partial_1, \partial_2$ denote partial derivatives with respect
to $\xi$, $\eta$, respectively.
\label{2:t1}
\end{thm}

We apply Theorem \ref{2:t1} to the solution
\[
u(\xi,\eta)= e^{-\xi^2/2} H_n(\xi)e^{\eta^2/2} H_n(i\eta),
\]
of (\ref{2:pde3}), and for $C$ we take the positively oriented boundary of the
rectangle $|\xi|\le \xi_1$, $|\eta|\le \eta_1$, where
$|\xi_0|<\xi_1$, $|\eta_0|<\eta_1$, so $n^+=n^-=1$.
Then let $\xi_1\to \infty$ and note that the integrals over the
vertical sides converge to~$0$. The integrals over the horizontal sides of
the rectangle give the same contribution because the integrand changes sign
when $(\xi,\eta)$ is replaced by $(-\xi,-\eta)$.  Therefore, we obtain,
for $\eta=\eta_1>|\eta_0|$,
\begin{equation}
\fl 2\pi u(\xi_0,\eta_0)=\int_{-\infty}^\infty
(v(\xi,\eta,\xi_0,\eta_0)\partial_2 u(\xi,\eta)-u(\xi,\eta)\partial_2v(\xi,\eta,\xi_0,\eta_0))d\xi.
\label{2:int2}
\end{equation}
Set $f(\eta)=f_n(\eta,\xi_0,\eta_0)$ and $g(\eta)= e^{\eta^2/2}H_n(i\eta)$.
Then (\ref{2:int2}) can be written as
\begin{equation}
2\pi u(\xi_0,\eta_0)=g'(\eta)f(\eta)-g(\eta)f'(\eta) .
\label{2:int3}
\end{equation}
By differentiating both sides of (\ref{2:int3}) with respect to $\eta$ we
see that $f$ satisfies (\ref{2:ode2}).
Since $f(\eta)$ goes to $0$ as $\eta\to\infty$, it follows that
\[
f(\eta)= c e^{-\eta^2/2} H_{-n-1}(\eta),
\]
where $c$ is a constant.
Going back to (\ref{2:int3}), we find that
\[
2\pi u(\xi_0,\eta_0)= c W,
\]
where $W=i^n$ is the (constant) Wronskian of $e^{-\eta^2/2}H_{-n-1}(\eta)$ and $g(\eta)$.
Therefore, if $\eta>|\eta_0|$, we obtain
\[
f_n(\eta,\xi_0,\eta_0) =2\pi (-i)^ne^{-\xi_0^2/2}H_n(\xi_0)e^{\eta_0^2/2}
H_n(i\eta_0)e^{-\eta^2/2}H_{-n-1}(\eta) .
\]
Using the well-known series expansion in terms of Hermite functions (Lebedev (1972) \cite[Theorem 2, page 71]{Lebedev}),
we obtain the following result.
\begin{thm}\label{2:t2}
For $\xi,\eta,\xi_0,\eta_0\in\R$ with $\eta>|\eta_0|$, 
\[
\fl K_0(r(\xi,\eta,\xi_0,\eta_0))=\displaystyle{\sqrt{\pi}}
e^{(\eta_0^2-\xi_0^2-\eta^2-\xi^2)/2}
\sum_{n=0}^\infty \frac{(-i)^n}{2^{n-1}n!}H_n(\xi)H_{-n-1}(\eta)H_n(\xi_0)H_n(i\eta_0),
\]
where $r$ is defined by (\ref{2:r}).
\end{thm}
The special case $\xi_0=\eta_0=0$ of Theorem \ref{2:t2} can be found in \cite[Problem 7, page 298]{Lebedev}.
If we multiply each $\xi,\eta,\xi_0,\eta_0$ by $\sqrt{k}$ then we obtain an
expansion for $K_0(k r)$.  Inserting this expansion into
(\ref{1:secondexpansion}) yields our final result.

\begin{thm}\label{2:t3}
Let $\bfx$, $\bfx_0$ be points on $\R^3$ with parabolic cylinder
coordinates $(\xi,\eta,z)$ and $(\xi_0,\eta_0,z_0)$, respectively.
If $\eta_\lessgtr:={\min \atop \max}\{\eta,\eta_0\}$ then
\begin{eqnarray*}
&&\fl\frac{1}{\|\bfx-\bfx_0\|}=\frac{2}{\sqrt{\pi}}\int_0^\infty
e^{-k/2(\xi^2+\eta^2+\xi_0^2-\eta_0^2)}\cos k(z-z_0)\\[0.2cm]
&&\hspace{2.0cm}\times\sum_{n=0}^\infty\frac{(-i)^n}{2^{n-1}n!}
H_n(\sqrt{k}\xi) H_{-n-1}(\sqrt{k}\eta_>) H_n(\sqrt{k}\xi_0)
H_n(i\sqrt{k}\eta_<)\,dk,
\end{eqnarray*}
where $\bfx\ne\bfx_0$.
\end{thm}

If we reverse the order of the infinite series and the definite integral
in the above expression we obtain
\begin{eqnarray*}
&&\fl\frac{1}{\|\bfx-\bfx_0\|}=\frac{2}{\sqrt{\pi}}
\sum_{n=0}^\infty\frac{(-i)^n}{2^{n-1}n!} \int_0^\infty
e^{-k/2(\xi^2+\eta^2+\xi_0^2-\eta_0^2)}
\\[0.2cm]
&&\hspace{2.0cm}
\times
H_n(\sqrt{k}\xi)
H_{-n-1}(\sqrt{k}\eta)
H_n(\sqrt{k}\xi_0)
H_n(i\sqrt{k}\eta_0)
\cos k(z-z_0) \,dk.
\end{eqnarray*}
It would be interesting to know the value of the above definite integral.

\section{Expansion of $J_0(kr)$ for parabolic cylinder coordinates}
\label{3}

Transforming equation (\ref{1:pde1}) to parabolic coordinates (\ref{2:pc})
we obtain
\[
\frac{\partial^2 u}{\partial \xi^2}+\frac{\partial^2 u}{\partial \eta^2}+k^2(\xi^2+\eta^2) u =0 .
\]
In order to simplify notation we will (temporarily) set $k=\frac12$ and $\zeta=i\eta$ with $\zeta$ real.
Thus we consider
\begin{equation}\label{3:pde2}
\frac{\partial^2 u}{\partial \xi^2}-\frac{\partial^2 u}
{\partial \zeta^2}+\frac14(\xi^2-\zeta^2) u =0 ,\quad \xi,\zeta\in\R.
\end{equation}
If $u_1(\xi)$ and $u_2(\zeta)$ are solutions of the
ordinary differential
equation (\ref{1:ode}) for some $\lambda$ then
$u(\xi,\zeta)=u_1(\xi)u_2(\zeta)$ solves (\ref{3:pde2}).

The function $W_{1/2}$ from (\ref{1:h1}) transformed to $(\xi,\zeta)$ becomes
\begin{equation}
w(\xi,\zeta,\xi_0,\zeta_0)=J_0\left(\frac12 \tilde r(\xi,\zeta,\xi_0,\zeta_0)\right),
\label{definitonw}
\end{equation}
where $\tilde r^2$ is a symmetric polynomial defined by
\begin{eqnarray*}
&&\fl 4\tilde r^2:= \left[(\xi-\xi_0)^2-(\zeta-\zeta_0)^2\right]
\left[(\xi+\xi_0)^2-(\zeta+\zeta_0)^2\right]\\[0.1cm]
&&\hspace{-1.83cm}=8\xi\xi_0\zeta\zeta_0+\xi^4+\xi_0^4+\zeta^4+\zeta_0^4 
-2\xi_0^2\zeta^2-2\xi_0^2\zeta_0^2-2\zeta^2\zeta_0^2-2\xi^2\xi_0^2-2\xi^2\zeta^2-2\xi^{2}\zeta_0^2,
\end{eqnarray*}
and $J_0$ is the order zero Bessel function of the first kind
(see (\ref{1:J})).
For fixed
$\zeta,\xi_0,\zeta_0\in\R$ consider the function $f:\R\to\R$ defined by
\begin{equation}
f(\xi):=w(\xi,\zeta,\xi_0,\zeta_0) .
\label{fdefinedbyw}
\end{equation}
We wish to expand this function in terms of (modified) parabolic cylinder harmonics
according to a general expansion theorem that is derived in the following subsection.

\subsection{Spectral theory of (modified) parabolic cylinder harmonics --
A singular Sturm-Liouville problem}
\label{31}

We discuss the Sturm-Liouville problem
\begin{equation}
-u''-\case14 x^2 u=\lambda u ,\quad -\infty<x <\infty,
\label{3:ode}
\end{equation}
involving the spectral parameter $\lambda$ subject to
\[ u\in L^2(-\infty,\infty) .\]
By replacing $\frac14 x^2$ by $-\frac14x^2$, we obtain the equation
describing the harmonic oscillator whose eigenfunctions (Hermite functions) and
the corresponding spectral theory is well-known. The
spectral problem associated with (\ref{3:ode}) is far less known.

A discussion of differential equation (\ref{3:ode}) and its solutions can be found in
\cite{Miller1955} and in section 8.2 of (\ref{3:ode})
Erd{\'e}lyi {\it et al.} (1982) \cite{ErdelyiHTFII}
(see also Magnus (1941) \cite{Magnus41};
Wells \& Spence (1945) \cite{WellsSpence};
Cherry (1948) \cite{Cherry48};
Darwin (1949) \cite{Darwin}).

We will follow Chapter 9 in Coddington \& Levinson (1955) \cite{CoddingtonLevinson}.
First note that, by \cite[Corollary 2, page 231]{CoddingtonLevinson},
equation (\ref{3:ode}) is in the limit-point case at $x=\pm\infty$.
Therefore, one can apply section 5 of Chapter 9 in \cite{CoddingtonLevinson}.

For $\lambda,x\in\C$ we define the functions
$u_1(\lambda,x)$, $u_2(\lambda,x)$ as the solutions of (\ref{3:ode}) uniquely determined
by the initial conditions
\[
u_1(\lambda,0)=u_2'(\lambda,0)=1,\quad u_1'(\lambda,0)=u_2(\lambda,0)=0.
\]
These functions may be expressed in terms of Kummer's function of the first kind
(\ref{Kummerfirst})
\begin{eqnarray}
u_1(\lambda,x)&=&e^{-\frac{i}{4}x^2}M\left(\case14+\case{i}{2}\lambda,
\case12,\case{i}{2}x^2\right),
\label{3:u1}
\\
u_2(\lambda,x)&=&e^{-\case{i}{4}x^2} xM\left(\case34+\case{i}{2}\lambda,
\case32,\case{i}{2}x^2\right).
\label{3:u2}
\end{eqnarray}
For $x>0$, the function
\begin{eqnarray}
\label{1:u3}
u_3(\lambda,x)&=&e^{-\case{i}{4}x^2} U\left(\case14+\case{i}{2}\lambda,\case12,
\case{i}{2}x^2\right)\\
&= &\frac{\sqrt\pi}{\Gamma\bigl(\case34+\case{i}{2}\lambda \bigr)} u_1(\lambda,x)
-(1+i)\frac{\sqrt\pi}{\Gamma\bigl(\case14+\case{i}{2}\lambda\bigr)} u_2(\lambda,x)\nonumber
\end{eqnarray}
is another solution of (\ref{3:ode}).
Here the Kummer function of the second
kind $U:\C\times\C\times(\C\setminus(-\infty,0])\to\C$
can be defined as (see Olver {\it et al.} (2010) \cite[(13.2.42)]{NIST})
\[
\fl U(a,b,z):=
\frac{\Gamma(1-b)}{\Gamma(a-b+1)}
M(a,b,z)
+
\frac{\Gamma(b-1)}{\Gamma(a)}
z^{1-b}
M(a-b+1,2-b,z).
\]
Except when $z=0$, each branch of $U$ is entire in $a$ and $b.$
We assume that $U(a,b,z)$ has its principal value.
The asymptotic behavior of the Kummer function of the second kind
\cite[(13.2.6)]{NIST} shows that $u_3(\lambda,\cdot)\in L^2(0,\infty)$
provided that $\Im\lambda<0$. Since (\ref{3:ode}) is in the limit-point case
at $+\infty$, $u_3$ is the only solution with this property except for a constant
factor.  The Titchmarsh-Weyl functions $m_{\pm\infty}(\lambda)$, $\Im\lambda\neq 0$,
are defined by the property that $u_1(\lambda,x)+m_{\pm\infty}(\lambda)u_2(\lambda,x)$
is square-integrable at $x=\pm\infty$.
Therefore, 
\[
m_\infty(\lambda)=
\left\{ \begin{array}{ll}
\displaystyle \displaystyle{ -(1+i)\frac{\Gamma(\frac34+\frac{i}{2}\lambda)}
{\Gamma(\frac14+\frac{i}{2}\lambda)}}& \quad\mathrm{if}\ \Im\lambda<0,\\[5pt]
\displaystyle \displaystyle{-(1-i)\frac{\Gamma(\frac34-\frac{i}{2}\lambda)}
{\Gamma(\frac14-\frac{i}{2}\lambda)}} & \quad\mathrm{if}\ \Im\lambda>0. \nonumber
\end{array} \right.
\]
By symmetry, we have
\[
m_{-\infty}(\lambda)=-m_{\infty}(\lambda).
\]
Using the notation of \cite[Theorem 5.1, page 251]{CoddingtonLevinson}, on obtains
\[
M_{11}(\lambda)=\frac{-1}{2 m_\infty(\lambda)},\quad M_{12}=M_{21}=0,
\quad M_{22}(\lambda)=\frac{m_{\infty}(\lambda)}{2}.
\]
Now using \cite[page 250, last line]{CoddingtonLevinson}, we find that, for $\lambda\in\R$,
\[
 \rho'_1(\lambda):=\rho'_{11}(\lambda)=\frac{1}{\pi} \lim_{\epsilon\to0+}
 \Im\left(M_{11}(\lambda+i\epsilon)\right)
=\frac{e^{\frac\pi2\lambda}}{4\sqrt2\pi^2}
\left|\Gamma(\case14+\case{i}{2}\lambda)\right|^2,
\]
since \cite[(5.4.5)]{NIST}
\[
\Gamma(\case{1}{4}+iy)\Gamma(\case34-iy)=\frac{\pi\sqrt2}{\cosh(\pi y)+i\sinh(\pi y)}.
\]
Moreover, $\rho_{12}(\lambda)=\rho_{21}(\lambda)=0$ and
\[
\rho'_2(\lambda):=\rho'_{22}(\lambda)=\frac{1}{\pi}
\lim_{\epsilon\to0+}\Im\left(M_{22}(\lambda+i\epsilon)\right)
=\frac{e^{\frac\pi2\lambda}}{2\sqrt2\pi^2}
\left|\Gamma(\case34+\case{i}{2}\lambda)\right|^2.
\]
Since the $\rho$-functions are real-analytic functions (with no jumps), we see
that the spectrum of (\ref{3:ode}) is the whole real line $\R$ and there are no
eigenvalues. The latter also follows from the known asymptotic behavior of the
solutions of (\ref{3:ode}) (see \cite[Chapter 12]{NIST}).

Using a variant of Stirling's formula (see (5.11.9) in Olver {\it et al.} (2010) \cite{NIST})
\begin{equation}\label{3:Stirling}
|\Gamma(x+iy)|\sim\sqrt{2\pi} |y|^{x-1/2}e^{-\pi|y|/2}\quad\textup{as $x,y\in\R$, $|y|\to\infty$},
\end{equation}
we can determine the
asymptotic behavior of $\rho'_j$, namely
\begin{eqnarray}
\label{form1}
\rho_1'(\lambda)&\sim& \frac{1}{2\pi}|\lambda|^{-1/2}e^{\pi(\lambda-|\lambda|)/2}\quad\textup{as $|\lambda|\to \infty$,}\\
\label{form2}
\rho_2'(\lambda)&\sim& \frac{1}{2\pi}|\lambda|^{1/2}e^{\pi(\lambda-|\lambda|)/2}\quad \textup{as $|\lambda|\to \infty$.}
\end{eqnarray}

Applying \cite[Theorem 5.2, page 251]{CoddingtonLevinson} to the
analysis above, we obtain the
following result on the spectral resolution associated with
equation (\ref{3:ode}).
\begin{thm}\label{3:t1}
For a given function $f\in L^2(\R),$ form the functions
\begin{equation}
g_j(\lambda)=\int_{-\infty}^\infty u_j(\lambda,x)
f(x)\, dx,\quad j=1,2,\ \lambda \in\R.
\label{gjlambda}
\end{equation}
Then
\[
g_j\in L^2(\R,\rho_j), \quad j=1,2,
\]
or, equivalently,
\[
\int_{-\infty}^\infty |g_j(\lambda)|^2\rho'_j(\lambda)\,d\lambda<\infty, \quad j=1,2.
\]
The function $f$ can be represented in the form
\begin{equation}
f(x) = \sum_{j=1}^2 \int_{-\infty}^\infty
u_j(\lambda,x)g_j(\lambda)\rho_j'(\lambda)\,d\lambda.
\label{3:rep}
\end{equation}
Moreover, we have Parseval's equation
\[
\int_{-\infty}^\infty |f(x)|^2\,dx= \sum_{j=1}^2 \int_{-\infty}^\infty
|g_j(\lambda)|^2\rho'_j(\lambda)\,d\lambda .
\]
\label{thmspectralparabolic}
\end{thm}

Equations (\ref{gjlambda}), (\ref{3:rep}) establish a one-to-one correspondence  between $f$ and $(g_1,g_2)$.
The integrals appearing in (\ref{gjlambda}),
(\ref{3:rep}) have to be interpreted in the $L^2$-sense. For instance, (\ref{gjlambda})
means that
\[
\int_{-n}^{n} u_j(\lambda,x)f(x)dx
\]
converges to $g_j(\lambda)$ in $L^2(\R,\rho_j)$ as $n\to\infty$.
Of course, if
\[
\int_{-\infty}^\infty \left|u_j(\lambda,x)f(x)\right|dx<\infty\quad \textup{for every $\lambda\in\R$,}
\]
then (\ref{gjlambda}) is also true pointwise.

\subsection{Applying the spectral theory to the expansion of $J_0(kr)$
in parabolic cylinder coordinates}

Since $f(\xi)=O(|\xi|^{-1})$ in (\ref{fdefinedbyw}) as $|\xi|\to\infty$,
we have that $f\in L^2(\R)$.
Therefore, we can expand $f$ using Theorem \ref{thmspectralparabolic},
according to (\ref{3:rep}).
For $\lambda\in\R$, we form the integrals
\begin{equation}
g_j(\lambda,\zeta, \xi_0,\zeta_0)=
\int_{-\infty}^\infty w(\xi,\zeta,\xi_0,\zeta_0)
u_j(\lambda,\xi)\,d\xi,\quad j=1,2 .
\label{2:g}
\end{equation}
These are absolutely convergent integrals because $f(\xi)=O(|\xi|^{-1})$
and $u_j(\lambda,\xi)=O(|\xi|^{-1/2})$ as $|\xi|\to\infty$.

Since $w(\xi,\zeta,\xi_0,\zeta_0)$ from (\ref{definitonw})
solves (\ref{3:pde2}) for fixed $(\xi_0,\zeta_0),$ and
$u_j(\lambda,\xi)$ from
(\ref{3:u1}), (\ref{3:u2})
solves (\ref{3:ode}), it follows from differentiation under the
integral sign followed by integration by parts
(see for instance Sch{\"a}fke (1963) \cite[Satz 8, page 26]{Schafke63})
that $g_j$ solves (\ref{3:ode}) as a function of $\zeta$.
Since the function $w$ is symmetric in its four variables, it is also true that $g_j$
solves (\ref{3:ode}) as function of $\xi_0$ for fixed $\zeta_0,\zeta$ and as function of
$\zeta_0$ for fixed $\xi_0,\zeta$.  From these properties of $g_j$ it follows easily that
there are functions
$c_{jk\ell m}:\R\to\R$ with $j,k,\ell,m=1,2$, depending
on $\lambda$ but not on $\zeta,\xi_0,\zeta_0$ such that
\begin{equation}
g_j(\lambda,\zeta,\xi_0,\zeta_0) = \sum_{k,\ell,m=1}^2
c_{jk\ell m}(\lambda)u_k(\lambda,\zeta)
u_\ell(\lambda,\xi_0)u_m(\lambda,\zeta_0).
\label{2:gg}
\end{equation}
This formula holds for all $\lambda,\zeta,\xi_0,\zeta_0\in\R$.
Substituting $\zeta=\xi_0=\zeta_0=0$ in (\ref{2:g}), (\ref{2:gg}), we obtain
\[
c_{j111}(\lambda)=\int_{-\infty}^\infty J_0(\case14 \xi^2)u_j(\lambda,\xi)\,d\xi.
\]
If $j=2,$ we integrate over an odd function, so $c_{2111}(\lambda)=0$.
By differentiating (\ref{2:gg}) with respect to $\zeta$ and/or $\xi_0$
and/or $\zeta_0$ and then substituting $\zeta=\xi_0=\zeta_0=0$ we find (after some calculations)
\[
c_{jk\ell m}(\lambda)=
\left\{ \begin{array}{ll}
\displaystyle c_1(\lambda) & \quad\mathrm{if}\ j=k=\ell=m=1,\\[2pt]
\displaystyle c_2(\lambda) & \quad\mathrm{if}\ j=k=\ell=m=2, \\[2pt]
\displaystyle 0 & \quad\textup{otherwise}. \nonumber
\end{array} \right.
\]
Here $c_j:\R\to\R$ for $j=1,2$ is given by (see \ref{5})
\begin{eqnarray}
&&\fl c_1(\lambda):=\int_{-\infty}^\infty J_0(\case14 \xi^2)
u_1(\lambda,\xi)\,d\xi
=\frac{2\sqrt{2}\pi e^{-\pi\lambda/2}}
{\cosh(\pi\lambda)|\Gamma\left(\frac34+\frac{i\lambda}{2}\right)|^2},
\label{3:c1}\\
&&\fl c_2(\lambda):=-\int_{-\infty}^\infty \xi^{-1}
J_1(\case14 \xi^2)u_2(\lambda,\xi)\,d\xi
=\frac{-4\sqrt{2}\pi e^{-\pi\lambda/2}}
{\cosh(\pi\lambda)|\Gamma\left(\frac14+\frac{i\lambda}{2}\right)|^2}.\label{3:c2}
\end{eqnarray}
According to (\ref{3:rep}), 
\begin{eqnarray}\label{3:wsumintegral}
\fl w(\xi,\zeta,\xi_0,\zeta_0)&=&\sum_{j=1}^2 \int_{-\infty}^\infty\!\!\!
c_j(\lambda)\rho'_j(\lambda) u_j(\lambda,\xi)
u_j(\lambda,\zeta) u_j(\lambda,\xi_0)u_j(\lambda,\zeta_0)d\lambda.
\end{eqnarray}
Note that
\begin{eqnarray}
c_1(\lambda)\rho_1'(\lambda)&=&\frac{1}{2\pi\cosh(\pi\lambda)} \left|\frac{\Gamma\left(\frac14+i\frac{\lambda}{2}\right)}{\Gamma\left(\frac34+i\frac{\lambda}{2}\right)}\right|^2
=\frac{1}{4\pi^3}\left|\Gamma\left(\case14+\case{i\lambda}{2}\right)\right|^4,
\label{3:crho1}\\
c_2(\lambda)\rho_2'(\lambda)&=& -\frac{2}{\pi\cosh(\pi\lambda)}\left|\frac{\Gamma\left(\frac34+i\frac{\lambda}{2}\right)}{\Gamma\left(\frac14+i\frac{\lambda}{2}\right)}\right|^2
=\frac{-1}{\pi^3}\left|\Gamma\left(\case34+\case{i\lambda}{2}\right)\right|^4,
\label{3:crho2}
\end{eqnarray}
where we used \cite[(5.4.4), (5.5.5)]{NIST}.
It follows from (\ref{3:Stirling}) that
\begin{eqnarray}
\label{c1lamr1lamp}
c_1(\lambda)\rho_1'(\lambda)&\sim & \frac{1}{\pi|\lambda|\cosh(\pi \lambda)}\quad\textup{as $|\lambda|\to\infty$,}\\
\label{c2lamr2lamp}
c_2(\lambda)\rho_2'(\lambda)&\sim& -\frac{|\lambda|}{\pi\cosh(\pi \lambda)}\quad\textup{as $|\lambda|\to\infty$}.
\end{eqnarray}
It is known (see Atkinson (1964), \cite[(8.2.5)]{Atkinson}) that, for fixed $x\in\C$, there is a constant $C$ such that $u_j(\lambda,x)=O(e^{C|\lambda|^{1/2}})$.
It follows from 
(\ref{c1lamr1lamp}),
(\ref{c2lamr2lamp}),
that the integrands in (\ref{3:wsumintegral}) decay exponentially and therefore 
the corresponding integrals are absolutely convergent.
By the identity theorem for analytic functions we see that equation (\ref{3:wsumintegral}) is true for all $\xi,\zeta,\xi_0,\zeta_0\in\C$.

After setting $\zeta=i\eta$ and $\zeta_0=i\eta_0$ in (\ref{3:wsumintegral}),
one obtains the following result.

\begin{thm}\label{3:t2}
Let $\xi,\eta,\xi_0,\eta_0\in\R$. Then
\begin{eqnarray*}
 && \hspace{-0.5cm}J_0\left(\frac12 r(\xi,\eta,\xi_0,\eta_0)\right)\\
   && \hspace{0.7cm}=\sum_{j=1}^2 \int_{-\infty}^\infty
c_j(\lambda)\rho'_j(\lambda) u_j(\lambda,\xi)
u_j(\lambda,i\eta) u_j(\lambda,\xi_0)u_j(\lambda,i\eta_0)\,d\lambda,
\end{eqnarray*}
where $r$ is given by (\ref{2:r}) and $c_j(\lambda)\rho_j'(\lambda)$ is 
given by (\ref{3:crho1}), (\ref{3:crho2}).
\end{thm}
In the special case $\xi_0=\eta_0=0$ (or correspondingly $\xi=\eta=0$), Theorem \ref{3:t2} can
be found in Buchholz (1953) \cite[(16), page 175]{Buchholz}.
Of course, if we multiply each $\xi,\eta,\xi_0,\eta_0$ by $\sqrt{2k}$ we get the expansion of
$J_0(k\, r(\xi,\eta,\xi_0,\eta_0))$.
This leads to the $J_0(kr)$ expansion of a fundamental solution
for the three-dimensional
Laplace equation in parabolic cylindrical coordinates.

\begin{thm}\label{3:t3}
Let $\bfx$, $\bfx_0$ be points on $\R^3$ with parabolic coordinates $(\xi,\eta,z)$ and $(\xi_0,\eta_0,z_0)$, respectively,
Then
\begin{eqnarray*}
&&\fl\frac{1}{\|\bfx-\bfx_0\|}=\sum_{j=1}^2 \int_0^\infty \int_{-\infty}^\infty
c_j(\lambda)\rho'_j(\lambda) \\[0.1cm]
&&\hspace{0.5cm}\times u_j(\lambda,2\sqrt{k}\xi)
u_j(\lambda,2i\sqrt{k}\eta) u_j(\lambda,2\sqrt{k}\xi_0)u_j(\lambda,2i\sqrt{k}\eta_0) e^{-k|z-z_0|}\,d\lambda\,dk .
\end{eqnarray*}
\end{thm}

\subsection{The Riemann method of integration}

The function $w(\xi,\zeta,\xi_0,\zeta_0)$ (as a function of $(\xi,\zeta)$) is a solution of the partial differential equation
(\ref{3:pde2}) and it satisfies the condition
$w(\xi,\zeta,\xi_0,\zeta_0)=1$ if $\xi-\xi_0=\pm(\zeta-\zeta_0)$.
This shows that $w$ is the Riemann function of
(\ref{3:pde2}); for instance, see Garabedian (1986) \cite{Garabedian}.

The Riemann method of integration applied to the partial differential
equation (\ref{3:pde2}) as in Volkmer (1980) \cite{Volkmer80}
gives, for all $\zeta,\xi_0,\zeta_0\in\R$,
\begin{eqnarray}
\label{3:eq1}
&&\fl 2u(\xi_0,\zeta_0)= u(\zeta-\zeta_0+\xi_0,\zeta)+u(-\zeta+\zeta_0+\xi_0,\zeta)\\
&&\int_{-\zeta+\zeta_0+\xi_0}^{\zeta-\zeta_0+\xi_0}
\left[w(\xi,\zeta,\xi_0,\zeta_0)\partial_2 u(\xi,\zeta)-\partial_2 w(\xi,\zeta,\xi_0,\zeta_0) u(\xi,\zeta)\right]\,d\xi ,\nonumber
\end{eqnarray}
where $u\in C^2(\R^2)$ is a solution of (\ref{3:pde2}). This formula 
for $\xi_0=\zeta=0$ and $u(\xi,\zeta)=u_1(\lambda,\xi)u_2(\lambda,\zeta)$ with $u_1,u_2$ from
(\ref{3:u1}), (\ref{3:u2}) (after replacing $\zeta_0$ by $\zeta$), implies that
\begin{equation}\label{3:eq2}
\int_{-\zeta}^\zeta J_0(\case14(\xi^2-\zeta^2))
u_1(\lambda,\xi)\,d\xi =2u_2(\lambda,\zeta).
\end{equation}
This equation allows us to transform the even solution $u_1$ into the odd
solution $u_2$ of equation (\ref{3:ode}).  One can prove (\ref{3:eq2}) directly by denoting the left-hand
side of (\ref{3:eq2}) by $f(\zeta)$ and then showing that
$f$ is an odd solution of (\ref{3:ode}) with $f'(0)=2$.

By differentiating (\ref{3:eq1}) first with respect to $\xi_0,\zeta_0$,
and using $w_1=u_2$, $w_2=u_1$, we obtain
(after a lengthy calculation) that
\[
\int_{-\zeta}^\zeta \frac{\xi\zeta}{\xi^2-\zeta^2}
J_1(\case14(\xi^2-\zeta^2))u_2(\lambda,\xi)\,d\xi =2u_1(\lambda,\zeta)-2u_2'(\lambda,\zeta).
\]
This formula can also be proved directly.

Let $h:\R^2\to\R$ be the function defined by
\[
h(\xi,\zeta):=
\left\{ \begin{array}{ll}
\displaystyle J_0\left(\case14(\xi^2-\zeta^2)\right) & \quad\mathrm{if}\ |\xi|<|\zeta|,\\[5pt]
\displaystyle 0 & \quad\textup{otherwise.} \nonumber
\end{array} \right.
\]
For fixed $\zeta$ this is an even function in $L^2(\R)$ which can be expanded according
to Theorem \ref{3:t1} (without knowing $c_j(\lambda)$),
so that
\[
\sign(\zeta)h(\xi,\zeta)=2\int_{-\infty}^\infty u_1(\lambda,\xi)u_2(\lambda,\zeta)\rho'_1(\lambda)\,d\lambda.
\]
If $\xi=0$, we obtain
\[
\sign(\zeta) J_0(\case14 \zeta^2)= 2\int_{-\infty}^\infty u_2(\lambda,\zeta)\rho_1'(\lambda)\, d\lambda.
\]
By Theorem \ref{3:t1}, this formulas allows us to conclude
\[
\int_0^\infty J_0(\case14\zeta^2)u_2(\lambda,\zeta)\,d\zeta
=\frac{\rho_1'(\lambda)}{\rho_2'(\lambda)},
\]
and this is in agreement with  (\ref{A:integral2}).

\section{Expansion of $J_0(kr)$ for elliptic cylinder coordinates}
\label{Jzeroellipticcylinder}

Consider equation (\ref{1:pde1}) for $k>0$ and elliptic coordinates on the plane
\begin{equation}\label{4:elliptic}
x= c\cosh\xi\cos\eta,\quad y=c \sinh\xi \sin\eta,
\end{equation}
where $\xi\in[0,\infty)$, $\eta\in\R,$ and $c>0$. Transforming $u(\xi,\eta)=U(x,y)$ we obtain
\begin{equation}\label{4:pde}
\frac{\partial^2 u}{\partial\xi^2}+\frac{\partial^2u}{\partial\eta^2} +k^2c^2(\cosh^2\xi-\cos^2\eta) u =0 .
\end{equation}
Separating variables $u(\xi,\eta)=u_1(\xi)u_2(\eta),$ leads to
\begin{eqnarray}
-u_1''(\xi)+(\lambda-2q\cosh 2\xi)u_1(\xi)&=&0 ,\label{4:modifiedMathieu} \\
\hspace{0.5cm}u_2''(\eta)+(\lambda-2q\cos 2\eta)u_2(\eta)&=&0 ,\label{4:Mathieu}
\end{eqnarray}
where $q=\frac14 c^2k^2>0$.

For $q\in\R$, Mathieu's equation (\ref{4:Mathieu}) is a Hill's differential equation with period $\pi$;
see \cite{NIST} or \cite{MeixnerSchafke54}. Here $q$ is positive
but in the next section $q$ will be negative.
As a Hill's equation, Mathieu's equation (\ref{4:Mathieu}) admits nontrivial $2\pi$-periodic solutions if and only if
$\lambda$ is equal to one of its eigenvalues $a_n(q)$, $n\in\N_0$ or $b_n(q)$, $n\in\N$.
If $\lambda=a_n(q),$ then Mathieu's equation has an even $2\pi$-periodic solution $\ce_n(\eta,q)$, and if
$\lambda=b_n(q)$ then Mathieu's equation has an odd $2\pi$-periodic solution $\se_n(\eta,q)$.
These functions are normalized according to
\[ \int_0^{2\pi} \ce_n^2(\eta,q)\,d\eta=\int_0^{2\pi} \se_n^2(\eta,q)\,d\eta =\pi .
\]
Moreover,
\[ \ce_n(\eta+\pi,q)=(-1)^n \ce (\eta,q),\quad  \se_n(\eta+\pi,q)=(-1)^n \se (\eta,q) .\]
Note that all solutions of Mathieu's equation are entire functions of $\eta$.

For these $2\pi$-periodic solutions of the Mathieu equation, we have the following expansion theorem.

\begin{thm}\label{4:t1}
Let $f(z)$ be a $2\pi$-periodic function that is analytic in an open doubly-infinite strip $S$ that contains the real axis.
Then
\begin{equation}\label{4:series}
 f(z)=\alpha_0\ce_0(z,q)+\sum_{n=1}^\infty(\alpha_n\ce_n(z,q)+\beta_n \se_n(z,q)),
\end{equation}
where
\[
\alpha_n= \frac1\pi\int_{-\pi}^{\pi} f(x) \ce_n(x,q)\,dx,\quad
\beta_n= \frac1\pi\int_{-\pi}^{\pi} f(x)\se_n(x,q)\,dx.
\]
The series (\ref{4:series}) converges absolutely and uniformly on any compact subset of the strip $S$.
\end{thm}

Let $(x_0,y_0)$ and $(x,y)$ be points on $\R^2$ with distance $r$.
Let $(x_0,y_0)$, $(x,y)$ have elliptic coordinates $(\xi_0,\eta_0)$ and $(\xi,\eta)$, respectively.
Then
\begin{eqnarray}\label{4:r}
\fl r^2=c^2\left[(\cosh\xi\cos\eta-\cosh\xi_0\cos\eta_0)^2
+(\sinh\xi\sin\eta-\sinh\xi_0\sin\eta_0)^2\right].
\end{eqnarray}
Clearly, $J_0(kr)$ as a function of $(\xi,\eta)$ solves (\ref{4:pde}).
We substitute $\zeta=i\xi$ and $\zeta_0=i\xi_0$. Then (\ref{4:pde}) changes to
\begin{equation}\label{4:pde3}
\frac{\partial^2 u}{\partial \zeta^2}-\frac{\partial^2 u}{\partial \eta^2}+ k^2c^2(\cos^2\eta-\cos^2\zeta) u =0,
\end{equation}
and $J_0(kr)$ is transformed to
\[
w(\zeta,\eta,\zeta_0,\eta_0)=J_0(k \tilde r),
\]
where
\begin{eqnarray*}
 \tilde r^2&=&c^2[(\cos\zeta\cos\eta-\cos\zeta_0\cos\eta_0)^2-(\sin\zeta\sin\eta-\sin\zeta_0\sin\eta_0)^2]\\
 &=& c^2(\cos(\zeta-\eta)-\cos(\zeta_0-\eta_0))(\cos(\zeta+\eta)-\cos(\zeta_0+\eta_0)) .
\end{eqnarray*}
The function $w(\zeta,\eta,\zeta_0,\eta_0)$ is an analytic function on $\C^4$ and it solves  equation (\ref{4:pde3}) as a function of $(\zeta,\eta)$.
It is the Riemann function of this differential
equation because $w(\zeta,\eta,\zeta_0,\eta_0)=1$ if $\zeta-\zeta_0=\pm(\eta-\eta_0)$.
For fixed $\eta, \zeta_0,\eta_0$ we wish to expand the function $\zeta\mapsto w(\zeta,\eta,\zeta_0,\eta_0)$ in a series of Mathieu functions according to
Theorem \ref{4:t1}.
To this end we have to evaluate the integral
\[
\int_{-\pi}^{\pi} w(\zeta,\eta,\zeta_0,\eta_0)\ce_n(\zeta,q)\,d\zeta,
\]
and a similar integral with $\ce_n$ replaced by $\se_n$.
Using Riemann's method of integration applied to a pentagonal curve, it has been shown in \cite{Volkmer83} that
\begin{eqnarray}
 \int_{-\pi}^\pi w(\zeta,\eta,\zeta_0,\eta_0)\ce_n(\zeta,q)\,d\zeta&=& \mu_n(q) \ce_n(\eta,q)\ce_n(\zeta_0,q)\ce_n(\eta_0,q), \label{4:int1}\\
 \int_{-\pi}^\pi w(\zeta,\eta,\zeta_0,\eta_0)\se_n(\zeta,q)\,d\zeta&=& \nu_n(q) \se_n(\eta,q)\se_n(\zeta_0,q)\se_n(\eta_0,q). \label{4:int2}
\end{eqnarray}
There do not exist explicit formulas for the quantities $\mu_n(q)$ and $\nu_n(q)$ but they can be determined as follows.
Mathieu's equation (\ref{4:Mathieu}) with $\lambda=a_n(q)$, $n\in\N_0$,
has the solution $u_1(\eta)=\ce_n(\eta,q)$. We choose a second linear independent solution $u_2(\eta)$.
Then there is $\sigma$ such that
\[ u_2(\eta+\pi)=\sigma u_1(\eta)+(-1)^n u_2(\eta) \]
and
\begin{equation}\label{4:mu}
\mu_n(q)=\frac{2(-1)^n \sigma}{W[u_1,u_2]} ,
\end{equation}
where $W[u_1,u_2]$ denotes the Wronskian of $u_1$ and $u_2$.
Similarly, Mathieu's equation (\ref{4:Mathieu}) with $\lambda=b_n(q)$, $n\in\N$,
has the solution $u_3(\eta)=\se_n(\eta,q)$. We choose a second linear independent solution $u_4(\eta)$.
Then there is $\tau$ such that
\[ u_4(\eta+\pi)=\tau u_3(\eta)+(-1)^n u_4(\eta) \]
and
\begin{equation}\label{4:nu}
\nu_n(q)=\frac{2(-1)^n \tau}{W[u_3,u_4]} .
\end{equation}
Now applying Theorem \ref{4:t1} and substituting $\zeta=i\xi$, $\zeta_0=i\xi_0,$ we obtain the following result.

\begin{thm}\label{4:t2}
Let $\xi,\eta,\xi_0,\eta_0 \in\C$, and let $k>0$, $c>0$, $q=\frac14c^2k^2$. Then
\begin{eqnarray*}
&&J_0(k r) = \frac1\pi \sum_{n=0}^\infty \mu_n(q)\ce_n(i\xi,q)\ce_n(\eta,q)\ce_n(i\xi_0,q)\ce_n(\eta_0,q) \\
&&\hspace{1.2cm}+\frac1\pi \sum_{n=1}^\infty
\nu_n(q)\se_n(i\xi,q)\se_n(\eta,q)\se_n(i\xi_0,q)\se_n(\eta_0,q) ,
\end{eqnarray*}
where $r$ is given by (\ref{4:r}).
\end{thm}
Theorem \ref{4:t2} agrees with expansion (23) (for $j=1$ and $\nu=0$), 
section 2.66 in Meixner \& Sch\"afke \cite{MeixnerSchafke54},
who have a slightly different notation. They use $\me_n(z,q)$, $n\in\Z$, where
\begin{eqnarray*}
\me_n(z,q)&:=&\sqrt{2}\ce_n(z,q) \quad \textup{if
$n\in\N_0,$}\\
\me_{-n}(z,q)&:=&-\sqrt2i\se_n(z,q)\quad  \textup{if $n\in\N.$}
\end{eqnarray*}
Moreover, the coefficients $\mu_n(q)$ and $\nu_n(q)$ are represented in a different form.
The proof of Theorem \ref{4:t2} based on Riemann's method of integration appears to be new.

We now use (\ref{1:firstexpansion}) to obtain our final result in this section.
\begin{thm}\label{4:t3}
Let $\bfx$, $\bfx_0$ be points on $\R^3$ with elliptic cylinder
coordinates $(\xi,\eta,z)$ and $(\xi_0,\eta_0,z_0)$, respectively.
Then
\begin{eqnarray*}
&&\fl\frac{1}{\|\bfx-\bfx_0\|}=
\frac1\pi \sum_{n=0}^\infty \int_0^\infty \mu_n(q)\ce_n(i\xi,q)\ce_n(\eta,q)\ce_n(i\xi_0,q)\ce_n(\eta_0,q)e^{-k|z-z_0|}\,dk\\
&&\hspace{-0.75cm}+\frac1\pi \sum_{n=1}^\infty \int_0^\infty
\nu_n(q)\se_n(i\xi,q)\se_n(\eta,q)\se_n(i\xi_0,q)\se_n(\eta_0,q)e^{-k|z-z_0|}\,
dk ,
\end{eqnarray*}
where $q=\frac14 c^2k^2$.
\end{thm}

\section{Expansion of $K_0(kr)$ for elliptic cylinder coordinates}
\label{ExpansionKzeroellipticcylinder}
Consider equation (\ref{1:pde2}) for $k>0$.
Transforming to elliptic coordinates (\ref{4:elliptic}), we obtain
\begin{equation}\label{5:pde}
\frac{\partial^2 u}{\partial\xi^2}+\frac{\partial^2u}{\partial\eta^2} -k^2c^2(\cosh^2\xi-\cos^2\eta) u =0 .
\end{equation}
Separating variables $u(\xi,\eta)=u_1(\xi)u_2(\eta),$ leads again to (\ref{4:modifiedMathieu}), (\ref{4:Mathieu})
but now $q=-\frac14 c^2k^2$ is negative.

We will need the following solutions of the modified Mathieu equation (\ref{4:modifiedMathieu}) when $q<0$; see \cite[\S 28.20]{NIST}.
Set $q=-h^2$ with $h>0$. For $n\in\N_0$,
$\Ie_n(\xi,h)$ is the even solution of (\ref{4:modifiedMathieu}) with $\lambda=a_n(q)$
with asymptotic behavior
\[ \Ie_n(\xi,h)\sim  I_n(2h\cosh \xi)\quad\textup{as $\xi\to+\infty$}, \]
while $\Ke_n(\xi,h)$ is the recessive solution determined by
\[  \Ke_n(\xi,h)\sim K_n(2h\cosh \xi)\quad\textup{as $\xi\to+\infty$}, \]
where $I_n(z):=i^{-n}J_n(iz)$ and $K_n(z)$ are the modified Bessel functions of the first
\cite[(10.27.6)]{NIST} and
second kinds respectively, with integer order $n$ (see (\ref{1:J}), (\ref{1:K})).
Similarly, for $n\in\N$, $\Io_n(\xi,h)$ is the odd
solution of (\ref{4:modifiedMathieu}) with $\lambda=b_n(q)$
with asymptotic behavior
\[ \Io_n(\xi,h)\sim  I_n(2h\cosh \xi)\quad\textup{as $\xi\to+\infty$}, \]
while $\Ko_n(\xi,h)$ is the recessive solution determined by
\[
\Ko_n(\xi,h)\sim K_n(2h\cosh \xi)\quad\textup{as $\xi\to+\infty$}.
\]

For fixed $\xi,\xi_0,\eta_0\in\R$ we wish to expand the
function
\[ v(\xi,\eta,\xi_0,\eta_0):=K_0(kr(\xi,\eta,\xi_0,\eta_0)) \]
with $r$ given by (\ref{4:r}) into a series of periodic Mathieu functions according to Theorem \ref{4:t1}.
The corresponding integrals appearing in the expansion will be computed based on the observation that
$(\xi,\eta)\mapsto v(\xi,\eta,\xi_0,\eta_0)$ is a fundamental solution of (\ref{5:pde}).
In fact, it is a solution of (\ref{5:pde}) and it has logarithmic singularities at the points $\pm (\xi_0,\eta_0+2m\pi)$,
where $m$ is any integer. Arguing as in Volkmer (1984) \cite[Theorem 1.11]{Volkmer84}, we have the following representation theorem for a solution of (\ref{5:pde}).

\begin{thm}\label{5:t1}
Let $u\in C^2(\R^2)$ be a solution of (\ref{5:pde}). Let $(\xi_0,\eta_0)\in\R^2$,
and let $C$ be a a closed rectifiable curve on $\R^2$ which does not pass through
any of the points $\pm (\xi_0,\eta_0+2m\pi)$, $m\in\Z$. Let $n^\pm_m$ be the winding number of $C$ with respect
to $\pm (\xi_0,\eta_0+2m\pi)$. Then we have
\begin{eqnarray*}
&&2\pi \sum_m \left[n^+_m u(\xi_0,\eta_0+2m\pi)+n^-_mu(-\xi_0,-\eta_0-2m\pi))\right]\\
&&\hspace{8mm}=\int_C (u\partial_2 v-v\partial_2 u)\,d\xi+(v\partial_1 u-u\partial_1 v)\,d\eta,
\end{eqnarray*}
where $\partial_1, \partial_2$ denote partial derivatives with respect
to $\xi$, $\eta$, respectively.
\end{thm}

In Theorem \ref{5:t1} we choose
\[ u(\xi,\eta)= u_1(\xi)u_2(\eta),\]
where
\[ u_1(\xi)=\Ke_n(\xi,h),\quad u_2(\eta)=\ce_n(\eta,q) .\]
Let $\xi_0>0$ and $\eta_0\in\R$. We take the curve $C$ to be the positively oriented boundary of the rectangle $\xi_1\le \xi\le \xi_2$, $\eta_0-\pi\le \eta\le \eta_0+\pi$, where
$|\xi_1|<\xi_0<\xi_2$. Consider the line integral $\int_C$ in Theorem \ref{5:t1}.
Since $u_2$ has period $2\pi,$ the line integrals along the horizontal segments of $C$ cancel each other.
When $\xi_2\to+\infty$ the asymptotic behavior of $u_1(\xi)$ shows that the 
integral along the right-hand vertical segment of $C$
tends to $0$ as $\xi_2\to+\infty$.
Therefore, setting
\[ f(\xi)=\int_{-\pi}^\pi v(\xi,\eta,\xi_0,\eta_0)u_2(\eta)\,d\eta=\int_{\eta_0-\pi}^{\eta_0+\pi} v(\xi,\eta,\xi_0,\eta_0)u_2(\eta)\,d\eta\]
for $|\xi|<\xi_0$, one obtains
\begin{equation}\label{5:eq1}
2\pi u_1(\xi_0)u_2(\eta_0)= u_1(\xi)f'(\xi)-u_1'(\xi)f(\xi).
\end{equation}
We now argue as in section \ref{2}.
By differentiating (\ref{5:eq1}) with respect to $\xi$, we find that $f(\xi)$ satisfies
the modified Mathieu equation (\ref{4:modifiedMathieu}). It is easy to see that $f(\xi)$ is an even function,
so $f(\xi)=c u_3(\xi)$, where $u_3(\xi)=\Ie_n(\xi,h)$ and $c$ is a  constant. Then (\ref{5:eq1}) implies that
\[
2\pi u_1(\xi_0)u_2(\eta_0) =c W[u_1,u_3].
\]
Since $W[u_1,u_3]=1$, 
\begin{equation}\label{5:eq2}
\frac1\pi\int_{-\pi}^\pi v(\xi,\eta,\xi_0,\eta_0)u_2(\eta)\,d\eta=2 u_1(\xi_0)u_2(\eta_0)u_3(\xi)\quad\textup{if $|\xi|<\xi_0$}.
\end{equation}
By the same reasoning, we see that (\ref{5:eq2}) is also true when $u_1(\xi)=\Ko_n(\xi,h)$, $u_2(\eta)=\se_n(\eta,q)$, $u_3(\xi)=\Io_n(\xi,h)$.

Expanding the function $\eta\mapsto v(\xi,\eta,\xi_0,\eta_0)$ according 
to Theorem \ref{4:t1}, the following result is obtained.
\begin{thm}\label{5:t2}
Let $\xi,\eta,\xi_0,\eta_0\in\R$ such that $|\xi|<\xi_0$, and let $k>0$, $c>0, q=-\frac14c^2k^2$, $h=\frac12ck$.
Then
\begin{eqnarray*}
 K_0(kr)&=&2 \sum_{n=0}^\infty \Ie_n(\xi,h)\ce_n(\eta,q)\Ke_n(\xi_0,h)\ce_n(\eta_0,q)\\
  &+&2 \sum_{n=1}^\infty \Io_n(\xi,h)\se_n(\eta,q)\Ko_n(\xi_0,h)\se_n(\eta_0,q),
\end{eqnarray*}
where $r$ is given by (\ref{4:r}).
\end{thm}

Theorem \ref{5:t2} agrees with Meixner \& Sch\"afke
(1954) \cite[section 2.66]{MeixnerSchafke54} although our notation and proof are different.

Inserting this result in (\ref{1:secondexpansion}), we obtain our final result.

\begin{thm}\label{5:t3}
Let $\bfx$, $\bfx_0$ be points on $\R^3$ with elliptic cylinder
coordinates $(\xi,\eta,z)$ and $(\xi_0,\eta_0,z_0)$, respectively.  If
$\xi_\lessgtr:={\min \atop \max}\{\xi,\xi_0\}$ then
\begin{eqnarray*}
&&\fl\frac{1}{\|\bfx-\bfx_0\|}=
\frac4\pi\sum_{n=0}^\infty \int_0^\infty \Ie_n(\xi_<,h)\ce_n(\eta,q)\Ke_n(\xi_>,h)\ce_n(\eta_0,q)\cos k(z-z_0)\,dk\\
&&\hspace{-0.75cm}+\frac4\pi\sum_{n=1}^\infty \int_0^\infty \Io_n(\xi_<,h)\se_n(\eta,q)\Ko_n(\xi_>,h)\se_n(\eta_0,q)\cos k(z-z_0)\,dk,
\end{eqnarray*}
where $q=-\frac14 c^2k^2$, $h=\frac12 ck$, and $\bfx\ne\bfx_0$.
\end{thm}

\appendix
\section{Integrals for the (modified) parabolic cylinder harmonics expansion
of $J_0(kr)$}
\label{5}

The following formulas are valid for $\Im \lambda<0$ :
\begin{eqnarray}
I_1:=\int_0^\infty J_0\left(\case14 \xi^2\right)
u_3(\lambda,\xi)\,d\xi&=& \frac12 \sqrt\pi(1-i)\frac{G_1}{G_2^2},
\label{A:int1}\\
I_2:=\int_0^\infty \xi^{-1}J_1\left(\case14 \xi^2\right)
u_3(\lambda,\xi)\,d\xi&=& -\sqrt\pi \left( \lambda\frac{1}{G_2}+2i\frac{G_2}{G_1^2}\right),
\label{A:int2}
\end{eqnarray}
where
\begin{eqnarray*}
G_1&=&G_1(\lambda)=\Gamma\left(\case14+\case{i}{2}\lambda\right),\\
G_2&=&G_2(\lambda)=\Gamma\left(\case34+\case{i}{2}\lambda\right).
\end{eqnarray*}
We believe that these integrals may be known but do not have a reference.
We will derive (\ref{A:int2}).
In (\ref{1:u3}) we use the integral representation \cite[(13.4.4)]{NIST}
\[
\Gamma(a)U(a,b,z)=\int_0^\infty e^{-zt}t^{a-1}(1+t)^{b-a-1}\,dt, \quad \Re z,\Re a>0.
\]
Substituting $4s=\xi^2$ and changing the order of integration, one obtains
\begin{equation}\label{A:eq1}
I_2=\frac{1}{2G_1} \int_0^\infty t^{-\frac34+\frac{i}{2}\lambda}
(1+t)^{-\frac34-\frac{i}{2}\lambda}
\int_0^\infty s^{-1}J_1(s)e^{-is(2t+1)}\,ds\,dt.
\end{equation}
From \cite[page 405]{Watson}, we have, for $t>0$,
\begin{eqnarray*}
\int_0^\infty s^{-1}J_1(s)\cos(s(2t+1))\,ds& =&0,\\
 \int_0^\infty s^{-1}J_1(s)\sin(s(2t+1))\,ds&=& t+(t+1)-2\sqrt{t}\sqrt{t+1} .
\end{eqnarray*}
Substituting these formulas in (\ref{A:eq1}), we can evaluate $I_2$ using three times
the formula for the beta function \cite[(5.12.3)]{NIST}
\[
B(z,w)=\frac{\Gamma(z)\Gamma(w)}{\Gamma(z+w)}=\int_0^\infty t^{z-1}(1+t)^{-z-w}\,dt ,\quad
\Re z,\Re w>0.
\]
This gives (\ref{A:int2}).

The proof of (\ref{A:int1}) is similar, but in (\ref{1:u3}) one should first use
\cite[(13.2.40)]{NIST}
\[
U(a,b,z)=z^{1-b} U(1+a-b,2-b,z).
\]

The formulas (\ref{A:int1}), (\ref{A:int2}) remain valid for real $\lambda$. By 
separating real and imaginary parts, we obtain for $\lambda\in\R$,
\begin{eqnarray}
\int_0^\infty J_0(\case14\xi^2)u_1(\lambda,\xi)\,d\xi &=&
\frac{\Re(G_1\overline{G}_2)+\Im(G_1\overline{G_2})}{|G_2|^2},\label{A:integral1}\\
\int_0^\infty J_0(\case14\xi^2)u_2(\lambda,\xi)\,d\xi &=&
\case12\left|\frac{G_1}{G_2}\right|^2,\label{A:integral2} \\
\int_0^\infty \xi^{-1}J_1(\case14\xi^2)u_1(\lambda,\xi)\,d\xi &=&
2\left|\frac{G_2}{G_1}\right|^2-\lambda ,\label{A:integral3}\\
\int_0^\infty \xi^{-1}J_1(\case14\xi^2)u_2(\lambda,\xi)\,d\xi &=&
2\frac{\Re(G_1\overline{G}_2)+\Im(G_1\overline{G_2})}{|G_1|^2}, \label{A:integral4}.
\end{eqnarray}
We may use
\[
\Re(G_1\overline{G}_2)+\Im(G_1\overline{G_2})
=\frac{\pi \sqrt2e^{-\frac12\pi\lambda}}{\cosh(\pi\lambda)}.
\]
Formulas (\ref{A:integral1}), (\ref{A:integral4}) give us the
integrals (\ref{3:c1}), (\ref{3:c2}) noting that we integrate
even functions in (\ref{3:c1}), (\ref{3:c2}).



\section*{Acknowledgements}
Part of this work was conducted while H.~S.~Cohl was a National Research Council
Research Postdoctoral Associate in the Information Technology Laboratory at the 
National Institute of Standards and Technology, Gaithersburg, Maryland, U.S.A.

\section*{References}

\end{document}